\documentclass[12pt]{article}
\usepackage[a4paper,left=1in,top=1in,width=6.3in,height=9in]{geometry}

\usepackage{amsmath,amsfonts,graphicx}
\usepackage[sort]{natbib}
\bibpunct{(}{)}{;}{a}{,}{,}

\numberwithin{equation}{section}

\newcommand{\tr}{\text{tr}}
\newcommand{\etr}{\text{etr}}
\newcommand{\diag}{\text{diag}}
\newcommand{\F}{\text{F}}
\newcommand{\Bing}{\text{Bing}}
\newcommand{\ACG}{\text{ACG}}
\newcommand{\FB}{\text{FB}}
\newcommand{\VM}{\text{VM}}
\newcommand{\WC}{\text{WC}}
\newcommand{\Unif}{\text{Unif}}
\newcommand{\MFB}{\text{MFB}}
\newcommand{\MF}{\text{MF}}
\newcommand{\MB}{\text{MB}}
\newcommand{\MBb}{\text{MB-bal}}
\newcommand{\MACG}{\text{MACG}}
\newcommand{\SO}{\text{SO}}

\begin{document}

\title{A new method to simulate the {B}ingham and related distributions in
  directional data analysis with applications}

\author{John T. Kent, Asaad M. Ganeiber and Kanti V. Mardia\\
Department of Statistics, University of Leeds}
\date{} %\date{\today}
\maketitle
\abstract{A new acceptance-rejection method is proposed and
  investigated for the Bingham distribution on the sphere using the
  angular central Gaussian distribution as an envelope.  It is shown
  to have high efficiency and to be straightfoward to use.  The method
  can also be extended to Fisher and Fisher-Bingham distributions on
  spheres and related manifolds.}

{\it Some key words: acceptance-rejection, angular central Gaussian
  distribution, Fisher distribution, Grassmann manifold, Kent
  distribution, special orthogonal group, Stiefel manifold, von Mises
  distribution.}

\section{Introduction}
Directional data analysis is concerned with statistical analysis on
various non-Euclidean manifolds, starting with circle and the sphere,
and extending to related manifolds. Comprehensive monographs are
available for statistical analysis in this setting; see,
e.g., \citet{FLE87,Mardia-Jupp00, Chikuse03}.  However, the subject of
simulation has received much less coverage, with the key contributions
scattered through the literature.

The need for effective simulation methods has grown in recent years as
directional distributions have become components in more
sophisticated statistical models, which are studied using MCMC methods.
For example, \citet{Green-Mardia06} used the matrix Fisher
distribution for random $3 \times 3$ rotation matrices in a
Bayesian model to align two unlabelled configurations of points in
$\mathbb{R}^3$, with an application to a problem of protein alignment
in bioinformatics.  

In general there are suitable direct methods of simulation, especially
methods based acceptance rejection, for the simpler directional
models.  However, it is necessary to resort to cumbersome MCMC methods
for the more complicated distributions.  The purpose of this paper is
to extend availability of acceptance rejection methods to a wider
class of directional distributions.  The starting point is a new
acceptance rejection method for the Bingham distribution, which can then
be used as a building block in a wider range of applications.

The paper is organized as follows.  Following some background and
preparation in Section 2, the new acceptance rejection simulation
method for the Bingham distribution is proposed and analyzed in
Section 3.  Special cases and extensions are covered in Sections 4 and
5.  Finally Section 6 sets the results of this paper in context by
reviewing the literature and summarizing the best available methods in
different settings.

The unit sphere $S_p = \{x \in \mathbb{R}^{q}: \ x^Tx = 1\}, \ p \geq
1$, comprises the unit vectors in $\mathbb{R}^q$, where throughout the
paper $p$ and $q$ are related by $q=p+1$.  The surface area of $S_p$
is given by $\pi_q = 2 \pi^{q/2}/\Gamma(q/2)$ and the differential
element of surface area can be written as $[dx]$.  Thus the uniform
distribution on $S_p$ can be written as $\pi_q^{-1} [dx]$.  A more
explicit formula can be given using polar coordinates.  For example,
the circle $S_1$ can be parameterized by $\theta \in [0, 2 \pi)$ with
uniform measure $d \theta /(2\pi)$.  The sphere $S_2$ can be
parameterized by colatitude $\theta \in [0,\pi]$ and longitude
$[0,2\pi)$ with uniform measure
\begin{equation}
\label{eq:unif-sphere}
\sin \theta d \theta d\phi /(4\pi).
\end{equation}

Strictly speaking a probability density on a manifold is a density
with respect to an underlying measure.  In Euclidean space
$\mathbb{R}^p$ the underlying measure is usually taken to be Lebesgue
measure $dx$ without explicit comment. But on other manifolds more
care is needed.  This paper is concerned with spheres and related
compact manifolds for which there is a natural underlying uniform
measure with a finite total measure.  To avoid repeated occurences of
normalizing constants such as $\pi_q$ and differential elements such
as $[dx]$, all such probability densities will be expressed with
respect to the uniform distribution. Thus we will write the density
for the uniform distribution on $S_2$ as $f(x) = 1$ (with respect to
itself) rather than as $f(x) = 1/(4\pi)$ (with respect to $[dx]$) or
as $f(x) = \sin \theta$ (with respect to $d\theta d \phi$).

\section{Background}
Recall the acceptance-rejection method of simulation.  Consider two
densities,
\begin{equation}
\label{eq:fg}
f(x) = c_f f^*(x), \quad g(x) = c_g g^*(x)
\end{equation}
where $f^*$ and $g^*$ are known functions, but where the normalizing
constants may or may not have a known explicit form. Suppose it is
possible to simulate easily from $g$ and it is desired to simulate
observations from $f$. The key requirement is that there is a known
bound of the form
\begin{equation}
\label{eq:bound}
f^*(x) \leq M^* g^*(x) \text{  for all } x
\end{equation}
for some constant $M^*$.  The acceptance-rejection algorithm proceeds
as follows.

\begin{itemize}
\item[Step]1. Simulate $X \sim g$ independently of $W \sim \text{Unif}(0,1)$.
\item[Step]2. If $W <  f^*(X)/\{M^* g^*(X)\}$, then accept $X$.
\item[Step]3. Otherwise go back to step 1.
\end{itemize}

\noindent {\em Comments}
\begin{itemize}

\item[(a)] If we set $M=c_f M^*/c_g$, then (\ref{eq:bound}) can be
  expressed equivalently as $f(x) \leq M g(x)$ for all $x$.

\item[(b)] The bound $M$ satisfies $M \geq 1$.  The number of trials
  needed from $g$ is geometrically distributed with mean $M \geq 1$.
  The efficiency is defined by $1/M$.  For high efficiency the bound
  $M$ should be as close to 1 as possible.

\item[(c)] The algorithm can be used even if the normalizing constants
  do not have a known explicit form.  However, to compute the
  efficiency analytically, it is necessary to know the normalizing
  constants.

\item[(d)] Suppose the density $g(x) = g(x; b)$ depends on a parameter
  $b$ with corresponding bound $M^*(b)$ in (\ref{eq:bound}).  If the
  normalizing constant $c_g = c_g(b) $ has a known explicit form, then
  it is possible to maximize the efficiency with respect to $b$, even
  if $c_f$ does not have a known explicit form.
\end{itemize}

When developing acceptance-rejection simulation methods for
directional distributions, there are several issues to consider:
\begin{itemize}
\item the need for good efficiency for a wide range of
  concentration parameters for $f$, ranging from uniform to highly concentrated.
  In similar problems on $\mathbb{R}^p$, the task is simpler when
  distributions are closed under affine transformations; in such cases
  it is sufficient to consider just a single standardized form of the
  distribution for $f$.

\item the challenge in finding a tractable envelope distribution.

\item the presence of trigonometric factors in the base measure when
expressed in polar coordinates, such as in  (\ref{eq:unif-sphere}).
\end{itemize}

Next, using the concavity of the log function, we give a general
inequality which is be useful in the construction of
acceptance-rejection algorithms.  Consider the function of $u \geq 0$,
\begin{equation}
\label{eq:concave}
\phi(u) = \frac{q}{2} \log(1+2u/b) - u - \frac{q}{2} \log(1+2u_0/b) + u_0,
\end{equation}
where $q>0$ and $0<b<q$ are fixed constants and $u_0 = (q-b)/2$.  The
last two terms on the righthand side of (\ref{eq:concave}) are
constants, chosen so that $\phi(u_0) = 0$.  The value of $u_0$ is
chosen so that the function $\frac{q}{2} \log(1+2u/b)$ has slope 1 at
$u=u_0$; hence $\phi'(u_0)=0$.  Also note that $\phi''(u) < 0$ for $u
\geq 0$ so that $\phi(u)$ is a concave function.  Therefore, $\phi(u)
\leq 0$ for all $u \geq 0$.  After exponentiating, this inequality can
be re-arranged as
\begin{equation}
\label{eq:exp-bound}
e^{-u} \leq e^{-(q-b)/2} \left(\frac{q/b}{1+2u/b}\right)^{q/2}.
\end{equation}

To illustrate the usefulness of (\ref{eq:exp-bound}),
we start in Euclidean space $\mathbb{R}^p, \ p \geq 1$,and construct
an acceptance-rejection algorithm for the multivariate normal
distribution using a multivariate Cauchy envelope.  The multivariate
normal distribution $N_p(0, \Sigma)$ has density
\begin{equation}
f(x) = c_f f^*(x), \quad f^*(x) = \exp(-\frac{1}{2} x^T \Sigma^{-1} x),\quad
c_f = c_f(\Sigma) = |2 \pi \Sigma|^{-1/2},
\end{equation}
for $x \in \mathbb{R}^p$.  The multivariate Cauchy distribution   $C_p(0, \Psi)$ has density
\begin{equation}
g(x) = c_g g^*(x), \quad g^*(x) = (1+x^T \Psi^{-1} x)^{-q/2}, \quad
c_g = c_g(\Psi) = \frac{\Gamma(q/2)}{\pi^{q/2}} |\Psi|^{-1/2},
\end{equation}
where for convenience we have substituted $q=p+1$.

If we set $\Psi = b \Sigma$ so that the scatter matrix for the Cauchy
is a scalar multiple of the covariance matrix for the normal, and if we set
$u=\frac{1}{2} x^T \Sigma^{-1} x$, then the inequality (\ref{eq:exp-bound})
leads to a bound on the densities with 
$$
M(b) = 2^{-p/2} q^{q/2} e^{-q/2}  \pi^{1/2} b^{-1/2} e^{b/2}/\Gamma(q/2).
$$
Mininizing over $0<b<q$ yields the optimal parameter $b=1$ with optimal bound
\begin{equation}
\label{eq:bound-n-c}
%M=2^{-p/2}  q^{q/2} e^{-p/2} \pi^{1/2}/\Gamma(q/2).
M = M(1) = \sqrt{2\pi e} \left(\frac{q}{2e}\right)^{q/2}\big/\,\Gamma(q/2), 
\quad q=p+1.
\end{equation}
Table \ref{table:eucl-eff} gives collection of the efficiencies $1/M$ as a function of
dimension $p$.  For large $p$, $M \sim \sqrt{qe/2}$ by Stirling's
formula.

\begin{table}
\begin{center}
\begin{tabular}{lccccccccc}
$p$ & 1 & 2 & 3 & 4 & 5 & 10 & 50 & 100\\
eff.& 66\%& 52\%& 45\%& 40\%& 36\%& 26\%& 12\%& 9\%
\end{tabular} \end{center}
\caption{Efficiency of A/R simulation method for the multiviate normal distribution in $p$ dimensions, using a multivariate Cauchy envelope.}
\label{table:eucl-eff}
\end{table}

Note that the efficiency declines slowly with the dimension, but is
still high enough to be feasible even for dimension $p=100$.
Of course, this is just a toy example since there are better ways to
simulate the normal distribution.  However, it is important for the
next section, both as a motiviating example and as a limiting case.

\section{Simulating the Bingham distribution with an 
angular central  gaussian envelope}
In this section we describe the ``BACG'' acceptance rejection method
to simulate the Bingham distribution using the angular central
Gaussian distribution as an envelope.  As before, let $q=p+1$.

The Bingham distribution, $\Bing_p(A)$ on $S_p$, $p \geq 1$, where the parameter
matrix $A$ is $q \times q$ symmetric, has density
\begin{equation}
\label{eq:bingham}
f_\Bing(x) = c_\Bing f_\Bing^*(x), \quad f_\Bing^*(x) = \exp(-x^T A x).
\end{equation}
The normalizing constant $c_\Bing = c_\Bing(A)$ can be expressed as a
hypergerometric function of matrix argument
\citep[p. 182]{Mardia-Jupp00}, but is not sufficiently tractable to be
of interest here.  The use of a minus sign in the exponent is
unconventional but simplifies later formulae.  Since $A$ and $A + cI$
define the same distribution for any real constant $c$, we may assume
without loss of generality that the eigenvalues of $A$ satisfy
\begin{equation}
\label{eq:evals}
0 = \lambda_1 \leq \lambda_2 \leq \cdots \leq \lambda_q.
\end{equation}

The angular central Gaussian distribution, ACG($\Omega)$ on $S_p$, where the
parameter matrix $\Omega$ is $q \times q$ symmetric positive definite,
takes the form
\begin{equation}
\label{eq:acg}
f_\ACG(x) = c_\ACG f_\ACG^*(x), \quad f_\ACG^*(x) = 
\left(x^T \Omega x\right)^{-q/2},
\quad c_\ACG = |\Omega|^{1/2}.
\end{equation}
The angular central Gaussian distribution is simple to simulate.  If
$y \sim N_q(0, \Sigma)$, where $\Sigma$ is positive definite, then $x
= y/||y|| \sim \text{ACG}(\Omega)$ with $\Omega = \Sigma^{-1}$
\citep[e.g.,][p. 182]{Mardia-Jupp00}.

Setting $u = x^T A x$ in (\ref{eq:exp-bound}) and setting $\Omega =
\Omega(b) = I + 2A/b, \ b>0,$ yields the envelope inequality on the starred
densities
\begin{equation}
\label{eq:b-acg}
\begin{split}
f_\Bing^*(x) &= e^{-u}\\
       &\leq e^{-(q-b)/2} \left(\frac{q/b}{1+2x^T A x/b}\right)^{q/2}\\
       &= e^{-(q-b)/2} \left(\frac{q/b}{x^T \Omega x}\right)^{q/2}\\
       &= e^{-(q-b)/2} (q/b)^{q/2} f_\ACG^*(x),
\end{split} \end{equation}
using the constraint $x^Tx = 1$.
The corresponding bound $M(b)$ takes the form 
\begin{equation}
\label{eq:bound-b-acg}
M(b) = c_\Bing e^{-(q-b)/2} (q/b)^{q/2} |\Omega(b)|^{-1/2}.
\end{equation}
It can be checked that $\log M(b)$ is convex in $b$ with unique
minimizing value given by the solution of
\begin{equation}
\label{eq:b-find}
\sum_{i=1}^q \frac{1}{b+2 \lambda_i} = 1,
\end{equation}
where the lefthand side of (\ref{eq:b-find}) ranges between $\infty$ and 0 as $b$ ranges
between 0 and $\infty$.  Let $b_0$ denote the solution to
(\ref{eq:b-find}) and let $M(b_0)$ denote the optimal bound.

It does not seem possible to evaluate $M(b_0)$ in a more useful form
analytically, but it is possible to say what happens asymptotically.
Replace $A$ in (\ref{eq:bingham}) by $\beta A$ and think of $A$ as a
fixed matrix as $\beta>0$ gets large.  Provided the $p$ largest
eigenvalues of $A$ are strictly positive, the ACG distribution (restricted to a hemisphere about the mode) 
converges to a $p$-dimensional multivariate Cauchy distribution, the
Bingham distribution converges to a $p$-dimensional multivariate
normal distribution, $b_0$ converges to 1 and the bound $M(b_0)$
converges to the bound (\ref{eq:bound-n-c}) where the dimensions $p$
and $q=p+1$ have the same meanings in both sections.

Empirically, it has been noticed that the limiting case is the worst
possible case.  For smaller values of the concentration matrix $A$,
the efficiencies will be higher.  Table \ref{table:bing-eff}
illustrates the pattern for $p=2$, i.e. $q=3$. The efficiency is never
lower than 52\%, the value from Table \ref{table:eucl-eff} for $p=2$.
This limiting value is attained in the concentrated bipolar case (when
$\lambda_2 = \lambda_3$ is large).  The girdle case ($\lambda_2 = 0$)
has higher efficiencies.  Each entry in this table has been
constructed from one million simulations, so that the standard errors
are negligible.

\begin{table}
\begin{center}
\begin{tabular}{rrr}
$\lambda_2$ & $\lambda_3$ & Efficiency \\
0 & 0 & 100\%\\
0 & 10 & 84\%\\
10 & 10 & 58\%\\
0 & 100 & 80\%\\
100 & 100 & 53\%\\
\end{tabular} \end{center}
\caption{Efficiency of the BACG A/R simulation method on $S_p$ with 
  $A = \diag(0, \lambda_2, \lambda_3)$ for the Bingham distribution with an
  ACG envelope.}
\label{table:bing-eff}
\end{table}

\section{Manifolds and models in directional data analysis}
In order to prepare for special cases and extensions of the Bingham
distribution, it is helpful to give a brief survey of some of the
various manifolds and models used in directional data analysis.  For each of
these manifolds there is a unique invariant measure which can be used to
define a uniform distribution.
\subsection{The sphere $S_p$ revisited}
A general model on the sphere $S_p$ is the Fisher-Bingham distribution
with density
\begin{equation}
\label{eq:fb}
f_\FB^*(x) = \exp(\kappa x^T \mu_0 - x^T A x) 
\end{equation}
where $\kappa \geq 0$, $\mu_0 \in S_p$ and $A (q \times q)$ is
symmetric, without loss of generality with smallest eigenvalue equal
to 0.  If $A=0$ the model reduces to the von Mises ($p=1$), the Fisher
($p=2$), or the von Mises-Fisher (any $p \geq 1$) distribution.  If
$\kappa=0$, the model reduces to the Bingham distribution considered
in Section 3.

\subsection{Real projective space $\mathbb{R}P_p$} 
Real projective space is defined as the quotient space $\mathbb{R}P_p
= S_p/\{1,-1\}$ in which two antipodal points or ``directions'' $\pm
x$ are identified with one another to represent the same ``axis''.
Since the Bingham and ACG densities have the property of antipodal
symmetry, $f(x) = f(-x)$, they can also be viewed as densities on
$\mathbb{R}P_p$.

\subsection{Complex projective space $\mathbb{C}P_p$}
Another quotient space of the sphere is complex projective space,
$\mathbb{C}P_p= S_{2p+1}/S_1$.  To understand this space, suppose a
unit vector $x \in \mathbb{R}^{2q}$, $q=p+1$, is partitioned as $x^T =
(x_1^T, x_2^T)$ where $x_1$ and $x_2$ are $q$-dimensional.  The
information in $x$ can also be represented by a $q$-dimensional
complex vector $z = x_1 +ix_2$.  Then $\mathbb{C}P_p$ is obtained from
$S_{2p+1}$ by identifying the scalar multiples $e^{i\theta} z$ with
one another for all $\theta \in [0, 2\pi)$.

If the $2q \times 2q$ symmetric concentration matrix $A$ for a
$\Bing_{2p+1}$ distribution can be partitioned in the form
\begin{equation}
  \notag
  A = \begin{bmatrix} A_1 & -A_2 \\ A_2 & \phantom{-} A_1 \end{bmatrix},
\end{equation}
where $A_1$ is symmetric and $A_2$ is skew symmetric, then then the
quadratic form $-x^TA x$ in the exponent of the Bingham density can be expressed
in complex notation as $-z^*A_C z$ where $A_C = A_1 + i A_2$.  In
terms of $z$, the density possesses complex symmetry, $f(z) =
f(e^{i\theta}z)$ for all $\theta \in [0, 2\pi)$.  When expressed in
complex notation this distribution is known as the complex
Bingham distribution $\text{CB}_p(A_C)$; it can also viewed as a
distribution on $\mathbb{C}P_p$ \citep{kc94b}.

\subsection{The special orthogonal group $\SO(r)$}
The special orthogonal group $\SO(r)$ is the space of $r \times r$
rotation matrices, $SO(r) = \{X \in \mathbb{R}^{r \times r}: X^TX=I_r,
\ |X|=1\}$.  A natural parametric distribution is given by the matrix
Fisher distribution $MF_r(F)$, with $r \times r$ parameter matrix $F$.
The density is given by
\begin{equation}
\label{eq:MF-SO}
f^*(X) = \exp\{ \tr(F^TX)\}.
\end{equation}
To describe the concentration properties of this distribution, it is
helpful to give $F$ a signed singular value decomposition
\begin{equation}
\label{eq:ssvd}
 F = U \Delta V^T.
\end{equation}
The adjective ``signed'' means that $U$ and $V$ are $r \times r$
rotation matrices and the elements of the diagonal matrix $\Delta$
satisfy $\delta_1 \geq \cdots \geq \delta_{r-1} \geq |\delta_r|$,
where the final element is negative if and only if $|F|<0$.

\subsection{The Stiefel manifold $\mathcal{V}_{r,q}$}
Let $1 \leq r \leq q$ and define the the Stiefel manifold $V_{r,q}=
\{X_1 \in \mathbb{R}^{q \times r}: X_1^TX_1 = I_r\}$ to be the space
of $q \times r$ column orthonormal matrices $X_1$, say.

The matrix Fisher-Bingham distribution on $V_{r,q}$, denoted
$MFB(F_1,A,C)$, with parameter matrices $F_1 (q \times r)$, $A (q \times
q \text{ symmetric})$ and $C (q \times q \text{ symmetric})$, is
defined by the density
\begin{equation}
\label{eq:MFB-stiefel}
f_{\MFB}(X) \propto \etr(F_1^TX_1-C X_1^TAX_1).
\end{equation}
Special cases include the matrix Fisher distribution, denoted
$MF(F_1)$, with density
\begin{equation}
\label{eq:MF-stiefel}
f_{\MF}(X) \propto \etr(F_1^TX_1),
\end{equation}
the ``full'' matrix Bingham distribution, denoted
$MB(A,C)$, with density
\begin{equation}
\label{eq:MB-stiefel-full}
f_{\MB}(X) \propto \etr(-C X_1^TAX_1),
\end{equation}
and the ``balanced'' matrix Bingham distribution, denoted
$MB(A)$, with density
\begin{equation}
\label{eq:MB-stiefel-balanced}
f_{\MBb}(X) \propto \etr(-X_1^TAX_1).
\end{equation}
If $r=q$ the balanced matrix Bingham distribution reduces to the
uniform distribution.  The reason is that $X_1$ is an orthogonal matrix
in this case, so that $X_1X_1^T = I_q$ and $\tr(X_1^TAX_1) = \tr(X_1
X_1^TA) = \tr(A)$ is constant in (\ref{eq:MB-stiefel-balanced}).

If $r=q-1$, it is possible to extend $X_1$ uniquely by adding an extra
column to form a $q \times q$ rotation matrix $X$, say.  Hence
$V_{q-1,q}$ can be identified with $SO(q)$.  However, the version of
the matrix Fisher distribution in (\ref{eq:MF-SO}), with an $r \times
r$ parameter matrix $F$, is more general than that in
(\ref{eq:MF-stiefel}), with an $r \times (r-1)$ parameter matrix
$F_1$.

If $r=q$, then $V_{q,q}$ is the same as the orthogonal group $O(q)$,
which is twice the size of the special orthogonal group $SO(q)$, since
in $O(q)$ a matrix $X$ is allowed to have determinant $\pm 1$.
Although the two densities (\ref{eq:MF-SO}) and (\ref{eq:MF-stiefel})
formally look the same, they live on different spaces.

\subsection{The Grassmann manifold $\mathcal{G}_{r,q}$} 
Let $1 \leq r < q$.  The Grassmann manifold $\mathcal{G}_{r,q}$ is
defined to be the set of all $r$-dimensional subspaces of
$\mathbb{R}^q$.  It can be described as a quotient space of a Stiefel
manifold $\mathcal{G}_{r,q} = V_{r,q}/O(r)$, in which a $q \times r$
column orthormal matrix $X_1$ is identified with $X_1R$ for all $r
\times r$ orthogonal matrices $R$.  It should be noted that the
notation for this manifold is not standardized; some authors write
$\mathcal{G}_{r,q-r}$ instead of $\mathcal{G}_{r,q}$.

Since $\tr(X_1^TAX_1) = \tr(R^TX_1^TAX_1R)$, the balanced matrix
Bingham distribution (\ref{eq:MB-stiefel-balanced}) on the Stiefel
manifold $\mathcal{V}_{r,q}$ can also be viewed as a distribution on
the Grassmannian manifold $\mathcal{G}_{r,q}$.

For every $r$-dimensional subspace in $\mathbb{R}^q$, there is a
unique complementary $(q-r)$-dimensional subspace.  If $X_1$ and $X_2$
are column orthonormal matrices, whose columns are bases of these
subspaces, then $X=(X_1 \ X_2)$ is a $q \times q$ orthogonal matrix.
Further $X_1$ follows a balanced matrix Bingham distribution on
$\mathcal{G}_{r,q}$ with parameter matrix $A$ if and only if $X_2$
follows a balanced matrix Bingham distribution on
$\mathcal{G}_{q-r,q}$ with parameter matrix $-A$ (but be warned that
the eigenvalues of $-A$ will not have the standardized form in
(\ref{eq:evals})).  Hence for simulation purposes, we may without loss
of generality suppose that $r \leq q/2$.

\section{Accidental isomorphisms}
The two quotient manifolds $\mathbb{R}P_p$ and $\mathbb{C}P_p$ of direct
interest for the Bingham distribution due to the existence of
``accidental isomorphisms'' in which the quotient manifold becomes
identified with another familiar manifold through a quadratic mapping.
These isomorphisms are called ``accidental'' because there does not
seem to be any systematic pattern.  In each case the uniform
distribution on the quotient manifold maps to the uniform distribution on
the new manifold, and the Bingham distribution maps to a distribution
related to the von Mises-Fisher distribution on the new manifold.  The
implications for simulation are laid out in the next subsections.
\subsection{$\mathbb{R}P_1 = S_1$}
Euclidean coordinates on the circle can be represented in polar
coordinates by $x = (x_1, x_2)^T$ where $x_1 = \cos \theta, x_2 = \sin
\theta, \ \theta \in [0,2\pi)$.  Consider a two-to-one map to a new
circle defined by $\phi = 2\theta$, with Euclidean coordinates
$y=(y_1,y_2)^T$ where $y_1 = \cos \phi = x_1^2-x_2^2, \ y_2 = \sin
\phi = 2 x_1 x_2$.  Note that the antipodal directions $\theta, \
\theta+\pi$ map to the same value of $\phi$, so that the map is in
fact a one-to-one map between $\mathbb{R}P_1$ and $S_1$.  A quadratic
form in $x$ can be rewritten as
$$
x^T A x = \frac{1}{2}(a_{11}-a_{22})y_1 + a_{12}y_2 + \frac{1}{2}(a_{11}+a_{22}),
$$
which is a linear function of $y$.  Hence a Bingham distribution,
whose density is quadratic in $x$ on $\mathbb{R}P_1$ can be
identified with a von Mises distribution, whose density is linear in
$y$, on $S_1$.

Similarly, in the ACG density the quadratic form $x^T\Omega x$ becomes
a linear function of $y$, so the density in $y$ reduces to the wrapped
Cauchy density \citep[p. 52]{Mardia-Jupp00}.

Suppose $A$ is diagonal, $A=\diag(0, \lambda)$.  In
this case the dominant axis of the Bingham distribution is the
$x_1$-axis.  The corresponding von Mises density takes the form
$$
f_{\VM}(y) \propto \exp(\kappa y_1), \quad \kappa=\lambda/2,
$$
so that the corresponding von Mises density has its mode in the 
$y_1$-direction.  The corresponding wrapped Cauchy density, with 
$\Omega= I + 2A/b$, takes the form,
$$
f_{\WC}(y) = \frac{(1-\rho^2)}{1+\rho^2-2\rho y_1}
$$
where $\rho=(\beta-1)/(\beta+1)$ \citep[p. 52]{Mardia-Jupp00}.

Hence the simulation method for the Bingham distribution with an ACG
envelope can be recast as a simulation for the von Mises distribution
with a wrapped Cauchy envelope.  It turns out that this latter method
is identical to the proposal of \citet{Best-Fisher79}, even up to the
choice of the optimal tuning constant $b$.

\subsection{$\mathbb{C}P_1 = S_2$}
The complex projective space $\mathbb{C}P_{k-2}$ arises in the study
of shape for configurations of $k$ landmarks in the plane, and the
identification with $S_1$ when $k=3$ was used to visualize the shape space
for triangles of landmarks \citep{Kendall84}.  \citet{kc94b} showed
that the complex Bingham distribution on $\mathbb{C}P_1$ can be
identified with the Fisher distribution on $S_2$.

Motivated by this accidental isomorphism, \citet{kc06a} developed a
complex Bingham quartic (CBQ) distribution on $\mathbb{C}P_p, \ p \geq
1$.  When $p=1$, this distribution reduces to the $\FB_5$
distribution.  \citet{Ganeiber12} developed an effective and
reasonably efficient simulation method for the CBQ distribution for $p>1$.
However, since the technique is not based on an angular central
Gaussian envelope, details will not be given here.

\subsection{$\mathbb{R}P_3 = SO(3)$}

There is a quadratic mapping taking an unsigned 4-dimensional unit
vector $\pm x$ to a $3 \times 3$ rotation matrix $X=M(x)=M(-x)$, say.
More specifically
\begin{equation}
\label{S3toSO3}
M(x) = \begin{bmatrix}
x_1^2+x_2^2-x_3^3-x_4^2 & -2(x_1x_4-x_2x_3) & 2(x_1x_3+x_2x_4)\\
2(x_1x_4+x_2x_3) & x_1^2+x_3^2-x_2^2-x_4^2 & -2(x_1x_2-x_3x_4)\\
-2(x_1x_3-x_2x_4) & 2(x_1x_2+x_3x_4) & x_1^2+x_4^2-x_2^2-x_3^2
\end{bmatrix}
\end{equation}
\citep[p. 285]{Mardia-Jupp00}.  Further a random axis $\pm x$ on
$\mathbb{R}P_3$ follows a Bingham distribution if and only if the
corresponding random matrix $M(x)$ follows a matrix Fisher
distribution.  In particular, if $A=\Lambda$ is diagonal, then
$F=\Delta$ in (\ref{eq:ssvd}) is also diagonal with the parameters related by
\begin{equation}
\label{eq:B-MF}
\lambda_1=0, \quad 
\lambda_2 = 2(\delta_2+\delta_3), \quad
\lambda_3 = 2(\delta_1+\delta_3), \quad
\lambda_4 = 2(\delta_1+\delta_2).
\end{equation}
\citet{kb12b} gives some further details.

A simple way to simulate a rotation matrix from the matrix Fisher
distribution $MF_3(F)$ for a general parameter matrix $F$ with signed
singular value decomposition (\ref{eq:ssvd}) is given as follows.
Using the BACG method simulate $x$ from $\Bing_3(\Lambda)$ with
$\Lambda$ given by (\ref{eq:B-MF}), and let $M(x)$ denote the
corresponding rotation matrix using (\ref{S3toSO3}).  Then $U M(x)
V^T$ follows the matrix Fisher distribution $MF_3(F)$.  From Table
\ref{table:eucl-eff}, the efficiency will be at least 45\%.

\section{Simulating the Fisher-Bingham distribution 
with an angular  central gaussian envelope}
The von Mises-Fisher density on $S_p$ takes the form (\ref{eq:fb})
with $A=0$.  The elementary inequality $(1-y)^2 \geq 0$, with $y = x^T
\mu_0$ can be re-arranged to give
\begin{equation}
\label{eq:vmf-b}
\begin{split}
f_\F^*(x) &\leq \exp\left[\left(\kappa/2\right)\left\{\left(x^T \mu_0\right)^2 
+ 1\right\}\right]\\
&= \exp\left\{\kappa - \left(\kappa/2 \right)x^T A x \right\} \\
&= e^\kappa f_\Bing^*(x),
\end{split} \end{equation} where $A = I_q - \mu_0 \mu_0^T$.  Hence an
acceptance rejection simulation method for the von Mises-Fisher
distribution can be constructed using a Bingham envelope.

The two sides of (\ref{eq:vmf-b}) match when $x=\mu_0$ so that it is
not possible to get a tighter bound.  In relative terms, the two
starred densities are maximally different when $x=-\mu_0$.  This
difference matters most when $\kappa$ is large, when the efficiency of
acceptance-rejection with a Bingham envelope drops to $50\%$; the
efficiency rises to 1 as $\kappa \rightarrow 0$.  Empirically the
efficiency lies between these two extremes for intermediate values of
$\kappa$.

The inequality (\ref{eq:vmf-b}) can be combined with Section 3 to
provide a method to simulate the von Mises-Fisher distribution with an
ACG envelope.  Of course there is no need for a new method for the von
Mises-Fisher distribution.  Good methods are already available; see
the Section 8 for a discussion.  However, the bounds of this
section can be combined with the previous section to simulate the
Fisher-Bingham distribution with an ACG envelope.

The Fisher-Bingham distribution on $S_p$ takes the form (\ref{eq:fb})
and can be bounded by a Bingham density
$$
f_\FB^*(x) \leq
\exp(\kappa- x^T A^{(1)} x),
$$
where $A^{(1)} = A +(\kappa/2)(I-\mu_0\mu_0^T)$.  Then Section 3 
can be used to bound this Bingham density  by an ACG density.

Say the Fisher-Bingham distribution is ``aligned'' if $\mu_0$ is an
eigenvector of $A$ and if the density has its mode at $x=\mu_0$.  In
this case the Bingham envelope usually has an efficiency of at least
50\%, with the efficiency falling below this level only when the
density excessively flat at its mode and the concentration is high.
Under high concentration this situation corresponds to the case where
the limiting normal distribution would have a singular covariance
matrix.

\section{The balanced matrix Bingham distribtion}
The density for the balanced matrix Bingham distribution was given by
(\ref{eq:MB-stiefel-balanced}).  It can be viewed as either a density
on the Stiefel manifold $\mathcal{V}_{r,q}$ (which is invariant under
multiplication on the right by $r \times r$ orthogonal matrices, or on
the Grassmannian manifold $\mathcal{G}_{r,q}$.  The $q \times q$
concentration matrix $A$ has the same form as for the Bingham
distribution in Section 3.

The matrix ACG distribution, denoted $MACG_{r,q}(\Omega)$, where $\Omega$ is
a positive definite symmetric $q \times q$ matrix, is also lies on
$\mathcal{V}_{r,q}$.  It  is also invariant  under rotation on the right
and hence can also be viewed as a distribution on the Grassmann
manifold $\mathcal{G}_{r,q}$.  The density takes the form
$$
g_{MACG}^*(X) = |X^T \Omega X|^{-q/2}, \quad c_g = |\Omega|^{r/2}
$$
\citep[e.g.][p. 40]{Chikuse03}.  Simulations from this distribution
can be constructed as follows.  Let $Y$ be a $q \times r$ matrix whose
columns are independently normally distributed, $N_q(0, \Omega^{-1})$.
Set $X=X(X^TX)^{-1/2}$ using the symmetric square root of a positive
definite matrix.  Then $X \sim MACG_{r,q}(\Omega)$.

If $\Omega$ is related to $A$ by $\Omega=\Omega(b) = I_q + 2A/b$ as in
Section 3, then the balanced matrix Bingham density can be bounded by
the matrix ACG density by using the inequality in (\ref{eq:exp-bound})
$r$ times.  Namely, let the eigenvalues of $X^TAX$ be denoted $u_1
\geq \cdots \geq u_r \geq 0$.  Since $f_{\MBb}^*(X) = \exp(\sum u_i)$
and $f_{\MACG}^*(X) = \{\prod (1+2u_i/b)\}^{-q/2}$, applying
(\ref{eq:exp-bound}) $r$ times yields the envelope bound
\begin{equation}
\label{eq:bound-mb-macg}
M(b) = c_\MBb \{e^{-(q-b)/2} (q/b)^{q/2} |\Omega(b)|^{-1/2}\}^r.
\end{equation}
Optimizing over $b$ yields the same equation (\ref{eq:b-find}) as
before with the same value for the optimal value $b_0$.

The efficiency is expected to decline as $r$ increases.  However, as
noted before, we may restrict attention to the case $r \leq q/2$.
More numerical investigation is needed of the efficiency in this setting.

\section{Review and commentary on different simulation 
methods}
Since the simulation literature for directional distributions is
widely scattered, it is useful to summarize the best simulation
methods for various distributions of interest.  Table \ref{table:recc}
lists several common distributions on different spaces, together with the
recommended method of simulation.

Recently, some MCMC simulation methods on manifolds have been proposed
by \citet{Kume-Walker09} (Fisher-Bingham on $S_p$), \citet{Habeck09}
(matrix Fisher on SO(3)), \citet{Hoff09} (matrix Fisher-Bingham
distributions on Stiefel and Grassmann manifolds) and
\citet{Byrne-Girolami13} (more general distributions).  However, there
is still an ongoing investigation into the efficiency of different
MCMC methods, so that the table entry will just state ``MCMC'' when
there is not a simpler more specific recommendation.  Further details
are given in the following subsections.

\begin{table}[t!]
\label{table:recc}
\caption{Recommended simulation methods various distributions on different 
  directional spaces}
\begin{tabular}{lll}
Distribution & Space & Simulation method \vspace{.2cm} \\ 
von Mises-Fisher & $S_p$ & \citet{Wood87}\\
Bingham & $S_p$ or $\mathbb{R}P_p$ & BACG\\
complex Bingham & $S_{2p+1}$ or $\mathbb{C }P_p$ & \citet{kc04a} \\
complex Bingham quartic & $S_{2p+1}$ or $\mathbb{C }P_p$ & \citet{kc06a}\\
aligned Fisher-Bingham &  $S_p$ & BACG-based\\
general Fisher-Bingham &  $S_p$ & MCMC\\
matrix Fisher & $\mathcal{V}_{r,q}$ & MCMC\\
matrix Fisher & $SO(3)$ & BACG\\
matrix Fisher & $SO(p), \ p>3$ & MCMC\\
balanced matrix Bingham & $\mathcal{V}_{r,q}$ or  $\mathcal{G}_{r,q}$ & 
BACG\\
general matrix Bingham & $\mathcal{V}_{r,q}$ & MCMC\\
matrix Fisher-Bingham & $\mathcal{V}_{r,q}$ & MCMC\\
\end{tabular}
\end{table}

\subsection{Uniform distribution on $S_p$}
The simplest general method to simulate a uniform distribution on the
unit sphere $S_p$ is to set $x = u /||u||$ where $u \sim N_q(0,I_q)$, \ 
$q=p+1$.  In low dimensions there are sometimes simpler methods using
polar coordinates.  E.g. on the circle $S_1$, let $\theta \sim
\Unif(0,2\pi)$.  On the sphere $S_2$ with colatitude $\theta$ and
longitude $\phi$, let $\cos \theta \sim \Unif(-1,1)$ independently of
$\phi \sim \Unif(0,2\pi)$.

On the Stiefel manifold the easiest approach is to simulate $U_1 (q
\times r)$ with independent $N(0,1)$ entries, and set $X_1= U_1
(U_1^TU_1)^{-1/2}$ using the symmetric square root of a positive
definite matrix.

\subsection{von Mises-Fisher distribution $\F_p(\kappa, \mu_0)$ on
  $S_p$}
For general $p \geq 1$, the recommended method of simulation is an
acceptance/rejection method due to \citet{Ulrich84}, as modified by
\citet{Wood94}.  This method uses a fractional linear transformation
of a beta variate to provide an envelope for $u = x^T \mu_0$.  It
gives good efficiency across the whole range of values for $\kappa$.
In particular, for large $\kappa$ the distribution of $2(1-u)$ is
approximately the squared radial part of a multivariate normal distribution
under the von Mises-Fisher model and of a 
multivariate Cauchy distribution under the envelope model, mimicking
the efficiency calculations in (\ref{eq:bound-n-c}). 

Once the distribution of $u \in [0,1]$ has been simulated, it is
straightforward to whole von Mises-Fisher distribution by
incorporating a uniformly distributed random direction $y$, say, on
$S_{p-1}$ (so $y$ is a $p$-vector). More specifically, if $R = [R_1 \
\mu_0]$ is any $q \times q$ rotatation matrix whose last column equals
$\mu_0$, let $x = u \mu_0 + (1-u^2)^{1/2} R_1^T y$.

For $p=1$ the Ulrich-Wood method is essentially identical to the
\citet{Best-Fisher79} method.  One small exception to the
recommendation to use the Ulrich-Wood method is the case $p=2$ dimensions
when $u$ follows a truncated exponential distribution and can be
simulated more simply by the inverse method without any need for
rejection \citep[p. 59]{FLE87}.

\subsection{Bingham distribution $\Bing_p(A)$ on $S_p$ or 
$\mathbb{R}P_p$}
The BACG method developed in this paper is the first general-purpose
acceptance/rejection simulaion method for the Bingham distribution.
However, earlier methods have been discussed in the literature for
some special cases.  In particular if $p=1$, the BACG method reduces
to \citet{Best-Fisher79} method for the von Mises distribution as
noted in Section 6.

If $p=2$ and either $0=\lambda_1 <\lambda_2 = \lambda_3$ (bipolar
case) or $0=\lambda_1 = \lambda_2 < \lambda_3$ (girdle case), the
simulation problem can be reduced to a one-dimensional problem.
\citet{Best-Fisher86} developed effective envelopes in these cases,
with efficiencies broadly comparable to the BACG method.

If the eigenvalues appear in pairs then the methods for the complex
Bingham can be used.  \citet{kc04a} developed several simulation
methods that sometimes are better than BACG.

The BACG method here supersedes the MCMC method of \citet{Kume-Walker06}.

\subsection{Fisher-Bingham distribution $\FB(\kappa,\mu_0,A)$ 
on $S_p$}
The Kent (or $\FB_5$) distribution on $S_2$ is a special case of an
aligned Fisher-Bingham distribution.  An efficient simulation
algorithm for $\FB_5$ was developed by \citet{kb05a}.  Otherwise, for
other aligned $\FB$ distributions when $p=2$ or for other values of
$p$, the BACG-based method developed in Section 6 is the recommended
method. 

In particular, these methods supersede earlier acceptance-rejection
methods developed by \citet{Wood87} for various special
types of aligned Fisher-Bingham distribution on $S_2$.  They also
supersede the acceptance rejection method of \citet[Appendix
A4]{Scealy-Welsh11} for a higher-dimensional version of the Kent
distribution, for which the efficiency drops to 0 under high
concentration when $p>2$.  In addition they supersede the MCMC
method of \citet{Kume-Walker09} in the aligned case.

For non-aligned Fisher-Bingham distributions, it is difficult to make
any firm theoretical statements about the behvaviour of the algorithm
in Section 6.  However, under moderate concentration it is still
likely to be preferable to the MCMC methods of \citet{Kume-Walker09}.

\subsection {Matrix Fisher distribution $\MF(F)$ on $SO(r)$}
When $r=2$, $SO(2)$ is the same as $S_1$ and the matrix Fisher on
SO(2) is identical to the von Mises distribution on $S_1$, so no new
methodology is needed.

When $r=3$ the accidental isomorphism in Section 5.3 reduces this case
to the Bingham distribution on $S_2$, which can be simulated by the
BACG method.

Earlier methods to simulate the matrix Fisher distribution on SO(3),
now superseded by BACG, were based on MCMC algorithms.  These include
\citet{Green-Mardia06} and \citet{Habeck09}.

The cases $r>3$ are at least partly covered by the next subsection.

\subsection{Matrix Fisher-Bingham distribution $\MFB(F,A,C)$ on
$\mathcal{V}_{r,q}$}

For the general matrix Fisher-Bingham distribution $\MFB(F,A,C)$ on
$\mathcal{V}_{r,q}$, there is not yet a convenient and efficient A/R
algorithm other than for the balanced matrix Bingham case, where a
solution was given in Section 7.  However, the recent MCMC algorithms
of \citet{Hoff09} and \citet{Byrne-Girolami13} can deal this this
case.

\subsection{Product manifolds}
Finally, the main setting not covered in this paper is the class of
product manifolds where multivariate versions of directional models
can be defined.  There are a few special cases where acceptance
rejection methods are available (e.g. \citet[supplementary material]
{Mardia-etal06} for the
sine and cosine versions of a bivariate version of the von Mises
distribution on the torus), but in general MCMC methods are needed.
\bibliographystyle{dcu} \bibliography{simbingham}
\end{document}